\documentclass[11pt,letterpaper,reqno]{amsart}
\usepackage{tikz}
\usetikzlibrary{positioning, shapes.geometric, arrows.meta, calc}
\usepackage{amssymb}
\usepackage{amsmath}
\usepackage{amsthm}
\usepackage{amsfonts}
\usepackage{bbm}
\usepackage{enumitem}
\usepackage{xcolor}
\usepackage{pgfplots}
\pgfplotsset{compat=1.18} 
\usepackage{booktabs}
\usepackage{graphicx}
\usepackage[T1]{fontenc}
\usepackage{doi}
\usepackage{float} 
\usepackage{comment} 

\let\originalleft\left
\let\originalright\right
\renewcommand{\left}{\mathopen{}\mathclose\bgroup\originalleft}
\renewcommand{\right}{\aftergroup\egroup\originalright}

\newcommand\doublecheck{\textcolor{blue}{\checkmark\kern-0.5em\checkmark}}
\newcommand{\newvtheorem}[2]{\newtheorem{#1}[theorem]{\llap{\textnormal{\doublecheck}\ }#2}}

\theoremstyle{plain}
\newtheorem{theorem}{Theorem}[section]
\newvtheorem{vtheorem}{Theorem}

\newvtheorem{vlemma}{Lemma}
\newtheorem{conjecture}[theorem]{Conjecture}

\newvtheorem{vproposition}{Proposition}

\newvtheorem{vcorollary}{Corollary}

\newvtheorem{vclaim}{Claim}

\theoremstyle{definition}

\theoremstyle{remark}

\newvtheorem{vremark}{Remark}

\begin{document}

\title[A counterexample to a conjecture of S\'ark\"ozy]{A counterexample to a conjecture of S\'ark\"ozy on sums and products modulo a prime}

\author[Q.~Tang]{Quanyu Tang}
\address{School of Mathematics and Statistics, Xi'an Jiaotong University, Xi'an 710049, P.~R.~China}
\email{tang\_quanyu@163.com, tangquanyu@stu.xjtu.edu.cn}

\subjclass[2020]{11B75, 05C69}

\keywords{sum-product problem, finite fields, independent sets}

\begin{abstract}
Let $p$ be a prime and, for $A\subseteq \mathbb F_p$, define $A^\ast=(A+A)\cup(AA)$. S\'ark\"ozy conjectured that there exist constants $c>0$ and $p_0$ such that, for every prime $p>p_0$, every set $A\subseteq \mathbb F_p$ with $|A|>\left(\frac12-c\right)p$ satisfies $\mathbb F_p^\times\subseteq A^\ast$. We disprove this conjecture: for every odd prime $p\ge 5$, there exists a set $A\subseteq \mathbb F_p$ with $|A|=\frac{p-1}{2}$ such that $1\notin A^\ast$. Thus no positive constant $c$ can satisfy S\'ark\"ozy's conjecture. Conversely, if $|A|>\frac{p}{2}$, then $A+A=\mathbb F_p$. Therefore the sharp threshold is exactly $\frac12$.
\end{abstract}

\maketitle

\section{Introduction}\label{sec:Introduction}

Let $p$ be a prime and let $A\subseteq \mathbb F_p$. Following S\'ark\"ozy \cite{Sarkozy2001}, define
\[
A^\ast := (A+A)\cup(AA),
\]
where
\[
A+A=\{a+a':a,a'\in A\},
\qquad
AA=\{aa':a,a'\in A\}.
\]
Thus $A^\ast$ consists of those residue classes modulo $p$ that can be represented either as a sum of two elements of $A$ or as a product of two elements of $A$.

In his 2001 list of unsolved problems, S\'ark\"ozy \cite[Conjecture~65]{Sarkozy2001} stated the following conjecture.

\begin{conjecture}[S\'ark\"ozy \cite{Sarkozy2001}]\label{conj:65}
There exist constants $c>0$ and $p_0$ such that, if $p>p_0$ is prime and $A\subseteq \mathbb F_p$ satisfies \(|A|>\left(\frac12-c\right)p\), then \(\mathbb F_p^\times\subseteq A^\ast\). Moreover, S\'ark\"ozy asked whether for every $\varepsilon>0$ one may take $c=\frac14-\varepsilon$.
\end{conjecture}

S\'ark\"ozy also remarked that a simple example shows that one cannot take $c>1/4$. Problems of this type belong to the study of algebraic equations with restricted solution sets over finite fields. S\'ark\"ozy proved density results for the solvability of the equations
\[
a+b=cd,
\qquad
ab+1=cd
\]
over large subsets of $\mathbb F_p$, and Gyarmati and S\'ark\"ozy later developed this circle further for more general equations over finite fields; see \cite{Sarkozy2005,Sarkozy2008,GyarmatiSarkozyI,GyarmatiSarkozyII}. For a later density/Ramsey-type viewpoint, see \cite{CsikvariGyarmatiSarkozy2012}; see also \cite{Cilleruelo2012} for an elementary survey around related combinatorial questions in finite fields.

The purpose of this note is to show that the actual situation in Conjecture~\ref{conj:65} is much stronger in the negative direction: \emph{no} positive constant $c$ works. In fact, the exact threshold is $1/2$.

Our main result proves that, for every odd prime $p\ge 5$, there exists a set $A\subseteq \mathbb F_p$ with
\[
|A|=\frac{p-1}{2}
\]
such that even the single nonzero residue class $1$ does not belong to $A^\ast$. Hence Conjecture~\ref{conj:65} is false. On the other hand, the threshold $1/2$ is best possible, because if $|A|>p/2$, then in fact $A+A=\mathbb F_p$.

The construction is elementary. We encode the two forbidden relations
\[
a+a'=1
\qquad\text{and}\qquad
aa'=1
\]
by a graph on $\mathbb F_p$. The connected components of this graph can be described explicitly: there are two small exceptional components, possibly one additional $2$-vertex component arising from the roots of $x^2-x+1$, and all remaining components are $6$-cycles generated by the two involutions above. Choosing suitable independent sets in each component yields a set of cardinality $(p-1)/2$ avoiding both representations of $1$.

The six-term orbit appearing in our proof is classical: it is the cross-ratio (or anharmonic) orbit generated by the involutions
\[
x\mapsto 1-x,
\qquad
x\mapsto x^{-1}.
\]
See \cite{Limaye1972} for the classical projective-geometric background and \cite[\S2]{Lu2015} for a convenient finite-field/projective-line formulation of the same six-term orbit and its exceptional cases.

\section{Main result}

We begin with the easy upper bound for the threshold.\footnote{Throughout, the symbol \doublecheck{} indicates that the proof of the corresponding statement has been formalized in Lean~4.}

\begin{vproposition}\label{prop:half}
Let $p$ be a prime and let $A\subseteq \mathbb F_p$. If $|A|>p/2$, then $A+A=\mathbb F_p$. In particular, \(\mathbb F_p^\times\subseteq A^\ast\).
\end{vproposition}

\begin{proof}
Fix any $x\in \mathbb F_p$. The two sets
\[
A
\qquad\text{and}\qquad
x-A:=\{x-a:a\in A\}
\]
have the same cardinality $|A|>p/2$, so they cannot be disjoint. Hence there exist $a,a'\in A$ such that \(a=x-a'\). Therefore
\[
x=a+a'\in A+A.
\]
Since $x\in\mathbb F_p$ was arbitrary, we obtain $A+A=\mathbb F_p$.
\end{proof}

We now state and prove the main theorem of this note.

\begin{vtheorem}\label{thm:main}
For every odd prime $p\ge 5$, there exists a set $A\subseteq \mathbb F_p$ such that
\[
|A|=\frac{p-1}{2}
\qquad\text{and}\qquad
1\notin A+A,\quad 1\notin AA.
\]
\end{vtheorem}

\begin{proof}
Define a graph $G$ on the vertex set $\mathbb F_p$ as follows. Two distinct vertices $u,v$ are adjacent if
\[
u+v=1
\qquad\text{or}\qquad
uv=1,
\]
and a loop is placed at a vertex $u$ if
\[
2u=1
\qquad\text{or}\qquad
u^2=1.
\]
Thus a subset $A\subseteq \mathbb F_p$ satisfies
\[
1\notin A+A
\qquad\text{and}\qquad
1\notin AA
\]
if and only if $A$ is an independent set in $G$, where by an independent set we mean a set containing neither two adjacent vertices nor a looped vertex.

We classify the connected components of $G$. Observe first that every vertex $y\in \mathbb F_p$ has at most two neighbours in $G$: one coming from the equation $y+z=1$, namely $z=1-y$, and, when $y\neq 0$, one coming from the equation $yz=1$, namely $z=y^{-1}$. Moreover, a loop occurs precisely when $y=\frac12$ or $y=\pm1$. Thus, to verify that a given finite set is a connected component, it suffices to list the neighbours of each of its vertices. We write $a \sim b$ if $a$ and $b$ are adjacent.

\medskip
\noindent\textbf{1. Two exceptional components.}
The set $\{0,1\}$ is a connected component: there is an edge $0 \sim 1$ since $0+1=1$, and $1$ carries a loop because $1^2=1$.

Similarly, $\{2,\frac12,-1\}$ is a connected component: there is an edge
\[
2 \sim \frac12
\qquad\text{since}\qquad
2\cdot \frac12=1,
\]
and an edge
\[
2 \sim (-1)
\qquad\text{since}\qquad
2+(-1)=1.
\]
Moreover, $\frac12$ carries a loop because $2\cdot \frac12=1$, and $-1$ carries a loop because $(-1)^2=1$.

\medskip
\noindent\textbf{2. A possible $2$-vertex component.}
Consider
\[
f(X)=X^2-X+1.
\]
Let $\delta\in\{0,2\}$ be the number of roots of $f$ in $\mathbb F_p$; since the discriminant of $f$ is $-3\neq 0$ in $\mathbb F_p$ for $p\ge 5$, the polynomial $f$ has either $0$ or $2$ distinct roots. If $u$ is a root, then
\[
u^2-u+1=0
\quad\Longrightarrow\quad
u(1-u)=1,
\]
so \(u^{-1}=1-u\). Hence, if $u$ and $v$ are the two roots of $f$, then
\[
u^{-1}=1-u=v,
\qquad
v^{-1}=1-v=u,
\]
and therefore $\{u,v\}$ is a connected component of $G$. These roots are distinct from
\[
0,\ 1,\ -1,\ \frac12,\ 2,
\]
since for $p\ge 5$ one has
\[
f(0)=1,\quad
f(1)=1,\quad
f(-1)=3,\quad
f \left(\frac12\right)=\frac34,\quad
f(2)=3,
\]
all nonzero in $\mathbb F_p$.

\medskip
\noindent\textbf{3. All remaining components are $6$-cycles.}
Let
\[
x\in \mathbb F_p\setminus\left\{0,1,-1,\frac12,2\right\}
\]
and assume that \(x^2-x+1\neq 0\). Consider the rational transformations
\[
s(y)=1-y,
\qquad
t(y)=y^{-1}.
\]
A direct computation shows that
\[
s^2=t^2=(st)^3=\mathrm{id}
\]
as rational transformations. Hence every word in $s$ and $t$ reduces to one of
\[
\mathrm{id},\qquad s,\qquad t,\qquad st,\qquad ts,\qquad sts.
\]
It follows that every vertex reachable from $x$ by a path in $G$ lies in
\[
\Omega(x):=
\left\{
x,\ 1-x,\ \frac1x,\ 1-\frac1x,\ \frac1{1-x},\ \frac{x}{x-1}
\right\}.
\]
Conversely, these six vertices lie in the same connected component, since
\[
x \sim 1-x \sim \frac1{1-x} \sim \frac{x}{x-1} \sim 1-\frac1x \sim \frac1x \sim x.
\]
Here the first, third and fifth adjacencies come from the relation $u+v=1$, while the second, fourth and sixth come from $uv=1$.

We claim that the six vertices in $\Omega(x)$ are pairwise distinct. Suppose otherwise. Then there exist distinct
\[
g_1,g_2\in\{\mathrm{id},s,t,st,ts,sts\}
\]
such that \(g_1(x)=g_2(x)\). Then the nonidentity transformation
\[
g=g_2^{-1}g_1
\]
fixes $x$. The possible fixed-point equations are:
\[
s(y)=y
\iff
1-y=y
\iff
y=\frac12,
\]
\[
t(y)=y
\iff
y^{-1}=y
\iff
y^2=1
\iff
y=\pm 1,
\]
\[
st(y)=y
\iff
1-\frac1y=y
\iff
y^2-y+1=0,
\]
\[
ts(y)=y
\iff
\frac1{1-y}=y
\iff
y^2-y+1=0,
\]
\[
sts(y)=y
\iff
\frac{y}{y-1}=y
\iff
y(y-2)=0.
\]
Thus a nontrivial fixed point can only be one of
\[
0,\ 1,\ -1,\ \frac12,\ 2,
\]
or a root of $y^2-y+1=0$, contrary to the choice of $x$. Hence the six vertices in $\Omega(x)$ are pairwise distinct. Since every vertex reachable from $x$ lies in $\Omega(x)$ and every element of $\Omega(x)$ lies in the connected component of $x$, that component is exactly $\Omega(x)$. Finally, every $y\in \Omega(x)$ has precisely the two neighbours $1-y$ and $y^{-1}$ in $G$, and these are exactly the two neighbours displayed in the cycle above. Thus no extra edges occur, and the connected component of $x$ is a $6$-cycle.

\medskip
\noindent\textbf{4. Construction of a large independent set.}
We now choose an independent set $A$ in $G$ as follows:
\begin{itemize}[leftmargin=2em]
\item from the component $\{0,1\}$ choose the vertex $0$;
\item from the component $\left\{2,\frac12,-1\right\}$ choose the vertex $2$;
\item if the component $\{u,v\}$ exists, choose one of its two vertices;
\item from each $6$-cycle choose three alternating vertices.
\end{itemize}
This produces an independent set $A$. There are \(p-5-\delta\) vertices outside the components
\[
\{0,1\},\qquad \left\{2,\frac12,-1\right\},\qquad \{u,v\}\ \text{(if it exists)},
\]
so the number of $6$-cycles is
\[
\frac{p-5-\delta}{6}.
\]
Therefore
\[
|A|
=
1+1+\frac{\delta}{2}+3\cdot\frac{p-5-\delta}{6}
=
\frac{p-1}{2}.
\]
Since $A$ is independent in $G$, we have $1\notin A+A$ and $1\notin AA$.
\end{proof}

\begin{vcorollary}\label{cor}
Conjecture~\ref{conj:65} is false.
\end{vcorollary}

\begin{proof}
Suppose for contradiction that there exist constants $c>0$ and $p_0$ such that for every prime $p>p_0$, every set $A\subseteq \mathbb F_p$ with \(|A|>\left(\frac12-c\right)p\) satisfies \(\mathbb F_p^\times\subseteq A^\ast\). Choose a prime
\[
p>\max\left\{p_0,\,5,\,\frac1{2c}\right\}.
\]
For this prime, the construction in Theorem~\ref{thm:main} yields a set $A\subseteq \mathbb F_p$ with \(|A|=\frac{p-1}{2}\) and \(1\notin A^\ast\). But $cp>\frac12$, so \(\frac{p-1}{2}>\left(\frac12-c\right)p\). This contradicts the assumed property. Therefore no such constants $c$ and $p_0$ exist.    
\end{proof}
\begin{vremark}\label{rem}
The proof actually yields the exact extremal value
\[
\max\{|A|:A\subseteq\mathbb F_p,\ 1\notin A+A,\ 1\notin AA\}=\frac{p-1}{2}.
\]
Indeed, the upper bound follows already from Proposition~\ref{prop:half}, while Theorem~\ref{thm:main} gives the matching lower bound.
\end{vremark}

\section{Lean formalization}\label{sec:lean}

We used Aristotle\footnote{Aristotle is an AI system for automated theorem proving developed by Harmonic; see \cite{AristotlePaper}. For the project page and additional public information, see \cite{AristotleProject}.} to obtain formal proofs of the four principal mathematical statements in this paper. The formalization code is available at~\cite{Git}. The code was checked in the following environment:
\begin{itemize}
  \item Lean version: \texttt{leanprover/lean4:v4.28.0}.
  \item Mathlib commit: \texttt{8f9d9cff6bd728b17a24e163c9402775d9e6a365}.
\end{itemize}

\noindent
More precisely, the formalized file proves the following four statements corresponding to the main results of the paper:
\begin{itemize}
\item Proposition~\ref{prop:half} is formalized as \texttt{sumset\_eq\_univ}.
\item Theorem~\ref{thm:main} is formalized as \texttt{exists\_large\_avoiding\_set}.
\item Corollary~\ref{cor} is formalized as \texttt{sarkozy\_conjecture\_false}.
\item Remark~\ref{rem} is formalized as \texttt{exact\_extremal\_value}.
\end{itemize}

\noindent
For background on Lean and Mathlib, see~\cite{Mathlib,Lean4}.

\section*{Acknowledgments}

We thank Wouter van Doorn for his help with the formalization of the proofs in this paper. We also used Aristotle to assist in the preparation of this paper.

\end{document}